\documentclass[11pt,twoside]{article}
\usepackage{multicol}
\usepackage[greek,english]{babel}
\usepackage[vcentering,dvips]{geometry}
\usepackage{makeidx}
\usepackage{amsfonts}
\usepackage{amssymb}
\usepackage{eucal}
\usepackage{amsxtra}
\usepackage{color}
\usepackage{setspace}
\usepackage{times}
\usepackage{titlesec}
\usepackage{sectsty}
\usepackage{fancyhdr}
\usepackage{fancyheadings}
\allsectionsfont{\large}
\usepackage{url}
\usepackage{doi}
\usepackage{hyperref}
\hypersetup{
    colorlinks,%
    citecolor=black,%
    filecolor=black,%
    linkcolor=black,%
    urlcolor=black
}
\usepackage[T1]{fontenc}
\usepackage{selinput}
\SelectInputMappings{%
  adieresis={ä},
  eacute={é},
  Lcaron={Ľ},
}
\addtocounter{section}{0} \numberwithin{equation}{section}
\newtheorem{theorem}{Theorem}[section]
\newtheorem{proposition}[theorem]{Proposition}
\newtheorem{lemma}[theorem]{Lemma}
\newtheorem{remark}[theorem]{Remark}

\newcommand{\CP}[1][]{\ensuremath{{\mathbb{C}P^{n}}\; }}
\newcommand{\CH}[1][]{\ensuremath{{\mathbb{C}H^{n}}\; }}
\newcommand{\C}[1][]{\ensuremath{{\mathbb{C}^{n}}\;}}
\newcommand{\MN}[1][]{\ensuremath{{M_{n}(c)}\;}}
\newcommand{\MNc}[1][]{\ensuremath{{M_{2}(c)}\;}}
\newcommand{\CPc}[1][]{\ensuremath{{\mathbb{C}P^{2}}\; }}
\newcommand{\CHc}[1][]{\ensuremath{{\mathbb{C}H^{2}}\; }}
\newcommand{\rt}[1][]{\ensuremath{{S^{*}}\;}}

\newcommand{\M}{{\emph{M }}}
\newcommand{\K}{{\emph{P }}}

\hypersetup{urlcolor=blue }

\begin{document}
\title{\textbf{Conditions of parallelism of $^{*}$-Ricci tensor of three dimensional real hypersurfaces in non-flat complex space forms}}
\author{\textsc{Georgios Kaimakamis and Konstantina Panagiotidou}}
\date{}
\maketitle

\begin{abstract}
This paper focuses on the study of three dimensional real hypersurfaces in non-flat complex space forms whose $^{*}$-Ricci tensor satisfies conditions of parallelism. More precisely, extension of existing results concerning real hypersurfaces with vanishing, semi-parallel and pseudo-parallel $^{*}$-Ricci tensor in case of ambient space being the complex hyperbolic space are provided. Furthermore,  new results concerning $\xi$-parallelism of $^{*}$-Ricci tensor of real hypersurfaces in non-flat complex space forms are presented.
\end{abstract}


\rhead[\centering{G. Kaimakamis and K. Panagiotidou}]{\thepage}
\lhead[\thepage]{\centering{$^{*}$-Ricci tensor}}

\section{\textsc{Introduction}}
A \emph{complex space form} is an n-dimensional Kähler manifold of constant holomorphic sectional curvature \emph{c}. A complete and simply connected complex space form is complex analytically isometric to complex projective space $\mathbb{C}P^{n}$ if $c>0$, or to complex Euclidean space $\mathbb{C}^{n}$ if $c=0$ or to complex hyperbolic space $\mathbb{C}H^{n}$ if $c<0$. The complex projective and complex hyperbolic spaces are called \emph{non-flat complex space forms}, since $c\neq0$ and the symbol $M_{n}(c)$ is used to denote them when it is not necessary to distinguish them.

A real hypersurface \M is an immersed submanifold with real co-dimension one in $M_{n}(c)$. The Kähler structure ($J, G$), where $J$ is the complex structure and $G$ is the Kähler metric of $M_{n}(c)$, induces on \M an almost contact metric structure ($\varphi,\xi,\eta,g$). The vector field $\xi$ is called \emph{structure vector field} and when it is an eigenvector of the shape operator $A$ of \M the real hypersurface is called \emph{Hopf hypersurface} and the corresponding eigenvalue is $\alpha=g(A\xi,\xi)$.

The study of real hypersurfaces \M in $M_{n}(c)$ was initiated by Takagi, who classified homogeneous real hypersurfaces in $\mathbb{C}P^{n}$ and divided them into six types, namely ($A_{1}$), ($A_{2}$), ($B$), ($C$), ($D$) and ($E$) in \cite{T2}. These real hypersurfaces are Hopf ones with constant principal curvatures. In case of $\mathbb{C}H^{n}$ the study of real hypersurfaces with constant principal curvatures  was started by Montiel in \cite{Mo} and completed by Berndt in \cite{Ber}. They are divided into two types, namely ($A$) and ($B$), depending on the number of constant principal curvatures and they are homogeneous and Hopf hypersurfaces.

Many geometers have studied real hypersurfaces in non-flat complex space forms when certain geometric conditions are satisfied. An important condition is that of the shape operator $A$ commuting with the structure tensor $\varphi$. More precisely, the following Theorem owed to Okumura in case of $\mathbb{C}P^{n}$  (\cite{Ok}) and to Montiel and Romero in case of $\mathbb{C}H^{n}$ (\cite{MR}) plays an important role in the proof of other Theorems.

\begin{theorem}\label{Th-Ok-Mo-Ro}
Let M be a real hypersurface of $M_{n}(c)$, $n\geq2$. Then $A\varphi=\varphi A$, if and only if M is locally congruent to a homogeneous real hypersurface of type (A). More precisely\\
 in case of $\mathbb{C}P^{n}$\\
    $(A_{1})$   a geodesic hypersphere of radius r , where $0<r<\frac{\pi}{2}$,\\
    $(A_{2})$  a tube of radius r over a totally geodesic $\mathbb{C}P^{k}$,$(1\leq k\leq n-2)$, where $0<r<\frac{\pi}{2}.$\\
 In case of $\mathbb{C}H^{n}$\\
    $(A_{0})$   a horosphere in $ \mathbb{C}H^{n}$, i.e a Montiel tube,\\
    $(A_{1})$  a geodesic hypersphere or a tube over a totally geodesic complex hyperbolic hyperplane $\mathbb{C}H^{n - 1}$,\\
    $(A_{2}) $  a tube over a totally geodesic $\mathbb{C}H^{k}$ $(1\leq k\leq n-2)$.
\end{theorem}

Generally, the Ricci tensor \textit{S}, of a Riemannian manifold is given by the relation
\begin{equation*}
S(X,Y)=trace\{Z\rightarrow R(Z, X)Y\},
\end{equation*}
where $X$, $Y$ are tangent vectors on \M. The same definition holds for real hypersurfaces in non-flat complex space forms. Real hypersurfaces in $M_{n}(c)$, $n\geq2$, in terms of their Ricci tensor satisfying geometric conditions such as parallelism and commutativity with other tensor fields of real hypersurfaces have been studied. A review of known results concerning the Ricci tensor of the real hypersurfaces can be viewed in \cite{NR1}.

In \cite{Ham} Hamada, motivated by Tachibana`s work in \cite{Tach}, where the \textit{$^{*}$-Ricci tensor} of almost Hermitian manifolds is defined, introduced the latter notion in case of real hypersurfaces in non-flat complex space forms. Therefore, the $^{*}$-Ricci tensor $S^{*}$ is given by
\begin{equation}\label{*-Ricci tensor}
S^{*}(X,Y)=\frac{1}{2}trace(Z\rightarrow R(X, \varphi Y)\varphi Z),\nonumber\
\end{equation}
where $X$, $Y$ are tangent vectors on \M.

Motivated by the work that has been done in case of studying real hypersurfaces in terms of their Ricci tensor, the authors began to study real  hypersurfaces in non-flat complex space forms in terms of their $^{*}$-Ricci tensor. More precisely,
in \cite{A2} real hypersurfaces in \MNc with parallel $^{*}$-Ricci tensor, i.e. $(\nabla_{X}S^{*})Y=0$, for any tangent vectors $X$, $Y$ on \M were classified.  In \cite{KP2014} conditions of semi-parallel $^{*}$-Ricci tensor, i.e. $(R(X,Y)\cdot S^{*})Z=0$, and pseudo-parallel $^{*}$-Ricci tensor, i.e. $(R(X,Y)\cdot S^{*})Z=L\{[(X\wedge Y)\cdot S^{*}]Z\}$, with $L$ being a non-zero function, were studied for real hypersurfaces in \CPc.

The aim of the present paper is to provide an analytic proof and extension of the existing results included in Theorems 2 and 3 in \cite{KP2014} in case of real hypersurfaces in \CHc. More precisely, the following results are proved

\begin{theorem}\label{tsemiparallel}
The only real hypersurface with semi-parallel $^{*}$-Ricci tensor is the geodesic hypersphere in \CHc with $coth(r)=2$.
\end{theorem}

\begin{theorem}\label{tpseudoparallel}
Every real hypersurface in \MNc with pseudo-parallel $^{*}$-Ricci tensor is a Hopf hypersurface. Furthemore, M is locally congruent to either a real hypersurface of type (A) or to a Hopf hypersurface satisfying relation $A\xi=0$, with L constant.
\end{theorem}

Furthermore, in this paper the condition of \emph{$\xi$-parallel} $^{*}$-Ricci tensor, i.e.
\begin{eqnarray}\label{xiparallel}
(\nabla_{\xi}S^{*})X=0,\mbox{for any tangent vector $X$ on \M,}
\end{eqnarray}
is studied and the following Theorem is proved

\begin{theorem}\label{txiparallel}
Every real hypersurface in \MNc with $\xi$-parallel $^{*}$-Ricci tensor is a Hopf hypersurface. Moreover, M is locally congruent to i) a real hypersurface of type ($A$) or ii) to a Hopf hypersurface with $A\xi=0$ or iii) to a Hopf hypersurface whose principal curvatures corresponding to the holomorphic distribution are constant in the direction of $\xi$.
\end{theorem}

This paper is organized as follows: In Section 2 basic relations and results about real hypersurfaces in \MNc are given. In Section 3 analytic proofs of Theorems \ref{tsemiparallel} and \ref{tpseudoparallel} are presented. Finally, in Section 4 proof of Theorem \ref{txiparallel} is provided.

\section{\textsc{Preliminaries}}
Throughout this paper all manifolds, vector fields etc are assumed to be of class $C^{\infty}$ and all manifolds are assumed to be connected. Furthermore, in case of \CPc we have $c=4$ and in case of  \CHc we have $c=-4$.

Let \M be a real hypersurface without boundary immersed in a non-flat complex space form $(M_{n}(c),G)$ with complex structure $J$ of constant holomorphic sectional curvature $c$. Let $N$ be a locally defined unit normal vector field on \M and $\xi=-JN$ be the structure vector field of \M. For any vector field $X$ tangent to \M relation
\[JX=\varphi X+\eta(X)N\]
holds, where $\varphi X$ and $\eta(X)N$ are respectively the tangential and the normal
component of $JX$. The Riemannian connections
$\overline{\nabla}$ in $M_{n}(c)$ and $\nabla$ in \M satisfy the relation
\[\overline{\nabla}_{X}Y=\nabla_{X}Y+g(AX,Y)N,\]
where $g$ is the Riemannian metric induced from the metric $G$ and for any vector fields $X$, $Y$ on \M.

The \emph{shape operator} $A$ of the real hypersurface $M$ in $M_{n}(c)$ with respect to $N$ is defined by
\[\overline{\nabla}_{X}N=-AX.\]
An almost contact metric structure $(\varphi,\xi,\eta, g)$ is induced on \M from $J$ of $M_{n}(c)$, where $\varphi$  is a tensor field of type (1,1) and is called \emph{structure tensor} and $\eta$ is an 1-form. The following relations hold
\begin{eqnarray}\label{eq-1}
&&g(\varphi X,Y)=G(JX,Y),\hspace{20pt}\eta(X)=g(X,\xi)=G(JX,N),\nonumber\\
&&\varphi^{2}X=-X+\eta(X)\xi,\hspace{20pt}
\eta\circ\varphi=0,\hspace{20pt} \varphi\xi=0,\hspace{20pt}
\eta(\xi)=1,\nonumber\\
&&g(\varphi X,\varphi
Y)=g(X,Y)-\eta(X)\eta(Y),\hspace{10pt}g(X,\varphi Y)=-g(\varphi
X,Y).\label{eq-2}\nonumber\
\end{eqnarray}
Moreover, $J$ being parallel implies $\overline{\nabla} J=0$ and this leads to
\begin{eqnarray}\label{eq-3}
\nabla_{X}\xi=\varphi
AX,\hspace{20pt}(\nabla_{X}\varphi)Y=\eta(Y)AX-g(AX,Y)\xi.\nonumber\
\end{eqnarray}
    The ambient space $M_{n}(c)$ is of constant holomorphic sectional
curvature $c$ and this results in Gauss and Codazzi equations are respectively given by
\begin{eqnarray}\label{eq-4}
R(X,Y)Z=\frac{c}{4}[g(Y,Z)X-g(X,Z)Y+g(\varphi Y ,Z)\varphi
X
\end{eqnarray}
$$-g(\varphi X,Z)\varphi Y-2g(\varphi X,Y)\varphi
Z]+g(AY,Z)AX-g(AX,Z)AY,$$
\begin{eqnarray}\label{eq-5}
\hspace{10pt}
(\nabla_{X}A)Y-(\nabla_{Y}A)X=\frac{c}{4}[\eta(X)\varphi
Y-\eta(Y)\varphi X-2g(\varphi X,Y)\xi],
\end{eqnarray}
where $R$ denotes the Riemannian curvature tensor on \M and $X$, $Y$, $Z$ are any vector fields on \M.

The tangent space $T_{P}M$ at every point $P$ $\in$ \M is decomposed as
\begin{eqnarray}
T_{P}M=span\{\xi\}\oplus \mathbb{D},\nonumber\
\end{eqnarray}
where $\mathbb{D}=\ker\eta=\{X\;\in\;T_{P}M:\eta(X)=0\}$ and is called (\emph{maximal}) \emph{holomorphic distribution} (\emph{if $n\geq3$}).
 Due to the above decomposition the vector field $A\xi$ can be written
 \begin{eqnarray}\label{eq-7}
 A\xi=\alpha\xi+\beta U,\nonumber\
 \end{eqnarray}
 where $\beta=|\varphi\nabla_{\xi}\xi|$ and
 $U=-\frac{1}{\beta}\varphi\nabla_{\xi}\xi\;\in\;\ker(\eta)$ is a unit vector field, provided
 that $\beta\neq0$.

Next, the following results concern any non-Hopf real hypersurface \M in \MNc with local orthonormal basis $\{U, \varphi U, \xi\}$ at a point \K of \M.
\begin{lemma}\label{basic-lemma}
Let M be a non-Hopf real hypersurface in \MNc. The following relations hold on M
\begin{eqnarray}\label{lemma-1}
&&AU=\gamma U+\delta\varphi U+\beta\xi,\hspace{20pt} A\varphi U=\delta U+\mu\varphi U,\hspace{20pt}A\xi=\alpha\xi+\beta U\label{eq-8}\\
&&\nabla_{U}\xi=-\delta U+\gamma\varphi U,\hspace{20pt}
\nabla_{\varphi U}\xi=-\mu U+\delta\varphi U,\hspace{20pt}
\nabla_{\xi}\xi=\beta\varphi U,\label{eq-9}\nonumber\\
&&\nabla_{U}U=\kappa_{1}\varphi U+\delta\xi,\hspace{20pt}
\nabla_{\varphi U}U=\kappa_{2}\varphi U+\mu\xi,\hspace{20pt}
\nabla_{\xi}U=\kappa_{3}\varphi U,\label{eq-10}\nonumber\\
&&\nabla_{U}\varphi U=-\kappa_{1}U-\gamma\xi,\hspace{5pt}
\nabla_{\varphi U}\varphi U=-\kappa_{2}U-\delta\xi,\hspace{5pt}
\nabla_{\xi}\varphi U=-\kappa_{3}U-\beta\xi,\label{eq-11}\nonumber\
\end{eqnarray}
where $\alpha, \beta, \gamma,\delta,\mu,\kappa_{1},\kappa_{2},\kappa_{3}$ are
smooth functions on M and $\beta\neq0$.
\end{lemma}

\begin{remark}
The proof of Lemma \ref{basic-lemma} is included in \cite{PX3}.
\end{remark}

The Codazzi equation (\ref{eq-5}) for $X$ $\in$ $\{U, \varphi U\}$ and $Y=\xi$ because of Lemma \ref{basic-lemma} implies
\begin{eqnarray}
\xi\delta&=&\alpha\gamma+\beta\kappa_{1}+\delta^{2}+\mu\kappa_{3}+\frac{c}{4}-\gamma\mu-\gamma\kappa_{3}-\beta^{2},\label{eq1}\\
(\varphi U)\alpha&=&\alpha\beta+\beta\kappa_{3}-3\beta\mu,\label{eq2}\\
(\varphi U)\beta&=&\alpha\gamma+\beta\kappa_{1}+2\delta^{2}+\frac{c}{2}-2\gamma\mu+\alpha\mu,\label{eq3}
\end{eqnarray}

and for $X=U$ and $Y=\varphi U$
\begin{eqnarray}\label{eq4}
U\delta-(\varphi U)\gamma&=&\mu\kappa_{1}-\kappa_{1}\gamma-\beta\gamma-2\delta\kappa_{2}-2\beta\mu.
\end{eqnarray}

Similar calculations to those of Theorem 2 in \cite{IvRy} imply that the $^{*}$-Ricci tensor of \M in \MNc since the ambient space is of constant holomorphic sectional curvature $c$ and $n=2$ is given by

\begin{equation}\label{starriccitensor}
S^{*}X= - [c\varphi^{2}X+(\varphi A)^{2}X],\;\;\mbox{for $X$ $\in$ $TM$}.
\end{equation}
If \M is a non-Hopf real hypersurface in \MNc and $\{U,\varphi U,\xi\}$ is a local orthonormal basis of it at some point $P$, the $^{*}$-Ricci tensor for $X$ $\in$ $\{U,\varphi U,\xi\}$ due to (\ref{eq-8}) and (\ref{starriccitensor}) takes the form
\begin{equation}\label{A1}
S^{*}{\xi}=\beta\mu U-\beta\delta\varphi U,\;\;S^{*}U=(c+\gamma\mu-\delta^{2})U\;\;\mbox{and}\;\;S^{*}\varphi U=(c+\gamma\mu-\delta^{2})\varphi U.
\end{equation}

Finally, the following Theorem which in case of \CP is owed to Maeda \cite{Maeda} and in case of \CH is owed to Montiel \cite{Mo} (also Corollary 2.3 in \cite{NR1}) is provided.
\begin{theorem}\label{Ma-Mo}
Let M be a Hopf hypersurface in $M_{n}(c)$, $n\geq2$. Then\\
i) $\alpha$ is constant.\\
ii) If $W$ is a vector field which belongs to $\mathbb{D}$ such that $AW=\lambda W$, then
\begin{eqnarray}\label{eq-A}
(\lambda-\frac{\alpha}{2})A\varphi W=(\frac{\lambda\alpha}{2}+\frac{c}{4})\varphi W.\nonumber\
\end{eqnarray}
iii) If the vector field $W$ satisfies $AW=\lambda W$ and $A\varphi W=\nu \varphi W$ then
\begin{eqnarray}\label{eq-B}
\lambda\nu=\frac{\alpha}{2}(\lambda+\nu)+\frac{c}{4}.
\end{eqnarray}
\end{theorem}

\begin{remark}\label{three}
In case of three dimensional Hopf hypersurfaces we can always consider a local orthonormal basis $\{W,\varphi W, \xi\}$ at some point $P$ $\in$ \M such that $AW=\lambda W$ and $A\varphi W=\nu\varphi W$. So relation (\ref{eq-B}) holds. Furthermore, the $^{*}$-Ricci tensor for $X$ $\in$ $\{W,\varphi W, \xi\}$ satisfies the relation
\begin{equation}\label{A2}
S^{*}{\xi}=0,\;\;\;\;S^{*}W=(c+\lambda\nu)W\;\;\mbox{and}\;\;S^{*}\varphi W=(c+\lambda\nu)\varphi W.
\end{equation}
\end{remark}

\section{\textsc{Proof of Theorems \ref{tsemiparallel} and \ref{tpseudoparallel}}}

Before proving Theorems \ref{tsemiparallel} and \ref{tpseudoparallel} the extension of Theorem 5 in \cite{KP2014} in case of real hypersurfaces in \CHc is given. More precisely, we obtain the following Theorem

\begin{theorem}\label{vanishing}
The only real hypersurface with vanishing $^{*}$-Ricci tensor is the geodesic hypersphere in \CHc with $coth(r)=2$.
\end{theorem}

In order to prove that every real hypersurface in \MNc with vanishing $^{*}$-Ricci tensor, i.e. $S^{*}X=0$, for any $X$ $\in$ $TM$ is a Hopf one, we follow the same steps as in the proof of Theorem 5 in \cite{KP2014}. The case of Hopf hypersurfaces in \CPc with vanishing $^{*}$-Ricci tensor is also included in the above proof. So it remains to examine the case of real hypersurfaces in \CHc in order to complete the proof of Theorem \ref{vanishing} of the present paper.

Since \M is a Hopf hypersurface in \MNc Theorem \ref{Ma-Mo} and remark \ref{three} hold. Since $S^{*}=0$ relation (\ref{A2}) implies that
\[c+\lambda\nu=0.\]
The above relation taking into account relation (\ref{eq-B}) yields that the real hypersurface has constant principal curvatures and this leads to the conclusion that a real hypersurface with vanishing $^{*}$-Ricci tensor is locally congruent to a real hypersurface of type ($A$) or type ($B$).

The following matrix includes the eigenvalues corresponding to three dimensional real hypersurfaces in \CHc according to \cite{Ber}. The type ($A_{1,1}$) refers to a geodesic hypersphere and the type ($A_{1,2}$) refers to a tube over a totally geodesic complex hyperbolic hyperplane $\mathbb{C}H^{1}$.
\\
\\
\begin{tabular}{l|*{5}{c}r}
Type             & $\alpha$ & $\lambda$ & $\nu$ & $m_{\alpha}$ & $m_{\lambda}$ & $m_{\nu}$ \\
\hline
($A_{0}$) & 2 & 1 & - & 1 & 2 & -  \\
($A_{1,1}$) &2$\coth(2r)$  &$\coth(r)$  &-  &1  & 2  &-   \\
($A_{1,2}$) &2$\coth(2r)$  &$\tanh(r)$  &-  &1  & 2  &-   \\
($B$) & 2$\tanh(2r)$ & $\tanh(r)$ & $\coth(r)$ & 1 & 1 & 1  \\
\end{tabular}
\\
\\

Substitution of the above eigenvalues in relation $c+\lambda\nu=0$ and because of $c=-4$ leads to the conclusion that only the eigenvalues of the geodesic hypersphere satisfies the latter. Furthermore, the radius $r$ of the geodesic hypersphere satisfies the relation $coth(r)=2$.

\subsection{Semi-parallel $^{*}$-Ricci tensor}
The $^{*}$-Ricci tensor is called semi-parallel when $(R(X,Y)\cdot S^{*})Z=0$, where $R$ is the Riemannian curvature and acts as derivation on $S^{*}$. More analytically, the above relation is written
\begin{eqnarray}\label{sp}
R(X,Y)S^{*}Z-S^{*}(R(X,Y)Z)=0\Rightarrow R(X,Y)S^{*}Z=S^{*}(R(X,Y)Z),
\end{eqnarray}
where $X$, $Y$ and $Z$ are any tangent vectors on \M.

Let $\mathcal{N}$ be the open subset of \M such that
\[\mathcal{N}=\{P\;\in\;\M:\beta\neq0\;\;\mbox{in a neighborhood of $P$}\}.\]

The inner product of relation (\ref{sp}) for $X=U$, $Y=\varphi U$ and $Z=U$ with $\varphi U$, due to (\ref{eq-4}) and (\ref{A1}) yields
\[\delta=0,\]
and relation (\ref{A1}) becomes
\begin{eqnarray}\label{A12*}
S^{*}\xi=\beta\mu U,\;\;\;\;S^{*}U=(c+\gamma\mu)U\;\;\mbox{and}\;\;S^{*}\varphi U=(c+\gamma\mu)\varphi U.
\end{eqnarray}

Furthermore, relation (\ref{sp}) for $X=\varphi U$, $Y=\xi$ and $Z=\varphi U$ due to (\ref{eq-4}) and (\ref{A12*}) implies
\[\mu(\frac{c}{4}+\alpha\mu)=0\;\;\mbox{and}\;\;(c+\gamma\mu)(\frac{c}{4}+\alpha\mu)=0.\]

Suppose that $\frac{c}{4}\neq\alpha\mu$ then the first of the above relations implies that $\mu=0$ and the second due to the latter results in $c=0$, which is a contradiction.

Therefore, on $\mathcal{N}$ relation $\frac{c}{4}+\alpha\mu=0$ holds. The inner product of relation (\ref{sp}) for $X=U$, $Y=\xi$ and $Z=U$ with $U$ because of (\ref{eq-4}) and (\ref{A12*}) yields
\[\mu(\frac{c}{4}+\alpha\gamma-\beta^{2})=0.\]

If $\frac{c}{4}+\alpha\gamma\neq\beta^{2}$ then we obtain $\mu=0$ and relation $\frac{c}{4}+\alpha\mu=0$ leads to $c=0$, which is a contradiction. So on $\mathcal{N}$ relation $\frac{c}{4}+\alpha\gamma=\beta^{2}$ holds.

The structure Jacobi operator $l=R_{\xi}$ of a real hypersurface in $M_{n}(c)$, $n\geq2$ is defined by
\[lX=R_{\xi}X=R(X,\xi)\xi.\]
In case of non-Hopf hypersurfaces \M in $M_{2}(c)$ taking into account relations
(\ref{eq-4}) and (\ref{eq-8}) the structure Jacobi operator is given by
\[lU=(\frac{c}{4}+\alpha\gamma-\beta^{2})U+\alpha\delta\varphi U,\;\;l\varphi U=\alpha\delta U+(\frac{c}{4}+\alpha\mu)\varphi U\;\;\mbox{and}\;\;l\xi=0.\]
Since $\delta=0$, $\frac{c}{4}+\alpha\mu=0$ and $\frac{c}{4}+\alpha\gamma=\beta^{2}$ we obtain
\[lU=l\varphi U=l\xi=0.\]

 It is known that there do not exist real hypersurfaces in \MN, $n\geq2$, with vanishing structure Jacobi operator (see Lemma 9 \cite{IRyan2009}). Thus, $\mathcal{N}$ is empty and the following Proposition is proved

\begin{proposition}
Every real hypersurface in \MNc whose $^{*}$-Ricci tensor is semi-parallel is a Hopf hypersurface.
\end{proposition}

Since \M is a Hopf hypersurface Theorem \ref{Ma-Mo} and remark \ref{three} hold. The case of Hopf hypersurfaces in \CPc with semi-parallel $^{*}$-Ricci tensor has been analytically studied in \cite{KP2014}. It remains the case of Hopf hypersurfaces in \CHc with $c=-4$. Relation (\ref{sp}) for $X=W$, $Y=\xi$ and $Z=W$  and for $X=\varphi W$, $Y=\xi$ and $Z=\varphi W$ because of relations (\ref{eq-4}) and (\ref{A2}) implies
\begin{eqnarray}\label{sp12}
(\lambda\nu-4)(\alpha\lambda-1)=0\;\;\mbox{and}\;\;(\lambda\nu-4)(\alpha\nu-1)=0.
\end{eqnarray}

Combination of the above relations implies that
\[\alpha(\lambda-\nu)(4-\lambda\nu)=0.\]

Suppose that $\alpha(\lambda-\nu)=0$ then we have two cases either $\alpha=0$ or $\lambda=\nu$. If $\alpha=0$ then relation (\ref{eq-B}) implies $\lambda\nu=-1$. Substitution of the latter relation in the first of (\ref{sp12}) leads to $-5=0$, which is a contradiction. If $\lambda=\nu$ then the shape operator $A$ commutes with the structure tensor $\varphi$ and because of Theorem  \ref{Th-Ok-Mo-Ro} \M is locally congruent to a real hypersurface of type ($A$). Moreover, the combination of relations (\ref{eq-B}) and the first of (\ref{sp12}) implies $\lambda^{2}(\lambda^{2}-4)=0$. Because of the matrix in section 3 we conclude that $\lambda^{2}=4$ and this occurs in case of geodesic hypersphere in \CHc.

Finally, if $\lambda\nu=4$ then relation (\ref{A2}) implies that the $^{*}$-Ricci tensor vanishes and owing to Theorem \ref{vanishing} we conclude that \M is a geodesic hypersphere and this completes the proof of Theorem \ref{tsemiparallel}

\subsection{Pseudo-parallel $^{*}$-Ricci tensor}
The $^{*}$-Ricci tensor is called pseudo-parallel when $(R(X,Y)\cdot S^{*})Z=L\{[(X\wedge Y)\cdot S^{*}]Z\}$, where $R$ is the Riemannian curvature and acts as derivation on $S^{*}$ and $L$ is a non-zero function. More analytically, the above relation is written
\begin{eqnarray}\label{pp}
R(X,Y)S^{*}Z-S^{*}(R(X,Y)Z)=L\{g(Y,S^{*}Z)X-g(X,S^{*}Z)Y-S^{*}[g(Y,Z)X-g(X,Z)Y]\},
\end{eqnarray}
where $X$, $Y$ and $Z$ are any tangent vectors on \M.

We consider $\mathcal{N}$ be the open subset of \M such that
\[\mathcal{N}=\{P\;\in\;\M:\beta\neq0\;\;\mbox{in a neighborhood of $P$}\}.\]

The inner product of relation (\ref{pp}) for $X=U$, $Y=\varphi U$ and $Z=U$ with $\varphi U$ because of (\ref{eq-4}) and (\ref{A1}) yields
\[\delta=0,\]
and relation (\ref{A1}) becomes
\begin{eqnarray}\label{A1*}
S^{*}\xi=\beta\mu U,\;\;\;\;S^{*}U=(c+\gamma\mu)U\;\;\mbox{and}\;\;S^{*}\varphi U=(c+\gamma\mu)\varphi U.
\end{eqnarray}

Relation (\ref{pp}) for $X=U$, $Y=\varphi U$ and $Z=\xi$ because of (\ref{eq-4}) and (\ref{A1*}) yields
\[\mu=0.\]
Moreover, relation (\ref{pp}) for $X=\varphi U$, $Y=\xi$ and $Z=\varphi U$ due to (\ref{eq-4}) and (\ref{A1*}) implies
\[\frac{c}{4}=L.\]

Relation (\ref{pp}) for $X=U$, $Y=\xi$ and $Z=U$ due to (\ref{eq-4}), (\ref{A1*}), $\mu=0$ and $\frac{c}{4}=L$ yields
\[\alpha\gamma=\beta^{2}.\]

On $\mathcal{N}$ relation (\ref{eq1}), (\ref{eq2}), (\ref{eq3}) and (\ref{eq4}) because of $\delta=\mu=0$  become
\begin{eqnarray}
\gamma\kappa_{3}&=&\beta\kappa_{1}+\frac{c}{4},\nonumber\\
(\varphi U)\alpha&=&\beta(\alpha+\kappa_{3}),\nonumber\\
(\varphi U)\beta&=&\beta^{2}+\beta\kappa_{1}+\frac{c}{2},\nonumber\\
(\varphi U)\gamma&=&\kappa_{1}\gamma+\beta\gamma.\nonumber\
\end{eqnarray}

Differentiation of $\alpha\gamma=\beta^{2}$ with respect to $\varphi U$ and taking into account all the above relations results in $c=0$ which is a contradiction.

Thus, $\mathcal{N}$ is empty and the following Proposition is proved

\begin{proposition}
Every real hypersurface in \MNc whose $^{*}$-Ricci tensor is pseudo-parallel is a Hopf hypersurface.
\end{proposition}

Since \M is a Hopf hypersurface, Theorem \ref{Ma-Mo} and remark \ref{three} hold. The case of Hopf hypersurfaces in \CPc with pseudo-parallel $^{*}$-Ricci tensor has been extensively studied in Theorem 3 in \cite{KP2014}. It remains the case of Hopf hypersurfaces in \CHc with $c=-4$. Relation (\ref{pp}) for $X=W$, $Y=\xi$ and $Z=W$  because of relations (\ref{eq-4}) and (\ref{A2}) implies
\[(\lambda\nu-4)(\alpha\lambda-1-L)=0.\]

Suppose that $\lambda\nu=4$ then relation (\ref{A2}) yields $S^{*}X=0$, for any vector field $X$ tangent to \M. The only real hypersurface with vanishing $^{*}$-Ricci tensor because of Theorem \ref{vanishing} is the geodesic hypersphere in \CHc with $\coth(r)=2$.

Next case $L=\alpha\lambda-1$ is examined. Relation (\ref{pp}) for $X=\varphi W$, $Y=\xi$ and $Z=\varphi W$ because of (\ref{eq-4}) and (\ref{A2}) implies
\[(\lambda\nu-4)(\alpha\nu-1-L)=0.\]

Suppose that $\lambda\nu=4$, then relation (\ref{A2}) implies that $S^{*}=0$ and due to Theorem \ref{vanishing} \M is geodesic hypersphere. Secondly, if $L=\alpha\nu-1$ combination of the latter relation with $L=\alpha\lambda-1$ results in
\[\alpha(\lambda-\nu)=0.\]

Thus, on \M either $\alpha=0$ or $\lambda=\nu$. If $\alpha=0$ then \M is locally congruent to a real hypersurface in \CHc with $A\xi=0$ (for the construction of these real hypersurfaces see \cite{IR2009}). If $\lambda=\nu$ it implies that the shape operator $A$ commutes with the structure tensor $\varphi$ and because of Theorem \ref{Th-Ok-Mo-Ro} it is concluded that \M is locally congruent to a real hypersurface of type ($A$) in \CHc.

Conversely, it is easily proved that the $^{*}$-Ricci tensor of the previous real hypersurfaces in \CHc have pseudo-parallel $^{*}$-Ricci tensor and that $L$ is constant given by $L=\alpha\lambda-1$ and this completes the proof of Theorem \ref{tpseudoparallel}.

\section{\textsc{Proof of Theorem \ref{txiparallel}}}

Let \M be a real hypersurface in \MNc whose *-Ricci tensor is $\xi$-parallel. More analytically, relation (\ref{xiparallel}) is written
\begin{eqnarray}\label{xipa}
\nabla_{\xi}(S^{*}X)=S^{*}(\nabla_{\xi}X),\;\;\mbox{for any $X$ $\in$ $TM$.}
\end{eqnarray}

Let $\mathcal{N}$ be the open subset of $M$ such that
\[\mathcal{N}=\{P\;\;\in\;\;M:\;\beta\neq0\;\;\mbox{in a neighborhood of $P$}\}.\]

On $\mathcal{N}$ the inner product of relation (\ref{xipa}) for $X=\xi$ with $\xi$ and $\varphi U$ because of (\ref{A1}) and relations of Lemma \ref{basic-lemma} implies respectively
\begin{eqnarray}\label{xi1}
\delta=0\;\;\mbox{and}\;\;\mu\kappa_{3}=c+\gamma\mu .
\end{eqnarray}

So relation (\ref{A1}) becomes
\begin{eqnarray}\label{A1n}
S^{*}{\xi}=\beta\mu U,\;\;S^{*}U=(c+\gamma\mu)U\;\;\mbox{and}\;\;S^{*}\varphi U=(c+\gamma\mu)\varphi U.
\end{eqnarray}

The inner product of relation (\ref{xipa}) for $X=\varphi U$ with $U$ due to relation (\ref{A1n}) and relations of Lemma \ref{basic-lemma} yields
\[\mu=0.\]
Substitution of the above relation in the second of (\ref{xi1}) results in $c=0$ which is a contradiction. Therefore, the following Proposition has been proved.

\begin{proposition}
Every real hypersurface in \MNc with $\xi$-parallel $^{*}$-Ricci tensor is a Hopf hypersurface.
\end{proposition}

Since \M is a Hopf hypersurface Theorem \ref{Ma-Mo} and remark \ref{three} hold. Relation (\ref{xipa}) for $X=W$ due to relation (\ref{A2}) and $\nabla_{\xi}W=\kappa W$, where $\kappa=g(\nabla_{\xi}W,W)$ and $g(\nabla_{\xi}W, W)=g(\nabla_{\xi}W,\xi)=0$ implies
\[\xi(\lambda\nu)=0.\]

Differentiating relation (\ref{eq-B}) with respect to $\xi$ and taking into account the fact that $\alpha$ is constant and the above relation we lead to
\[\alpha[\xi(\lambda+\nu)]=0.\]

Suppose that $\alpha\neq0$ then the above relation implies $(\xi\lambda)=-(\xi\nu)$. Substituting the last one in relation $\xi(\lambda\nu)=0$ we obtain
\[(\lambda-\nu)(\xi\lambda)=0.\]

If $(\xi\lambda)\neq0$ then $\lambda=\nu$ and this results in $A\varphi=\varphi A$. The last relation because of Theorem \ref{Th-Ok-Mo-Ro} implies that $M$ is locally congruent to a real hypersurface of type ($A$).

If  $(\xi\lambda)=0$ then also $(\xi\nu)=0$ and since $\lambda, \nu$ are the principal curvature corresponding to the holomorphic distribution we conclude that $M$ is locally congruent to a  Hopf hypersurface with constant $\lambda,\nu$ in direction of $\xi$.

The remaining case is $\alpha=0$ which implies that $M$ is a Hopf hypersurface with $A\xi=0$. More analytically, in case of \CPc, $M$ is locally congruent to a geodesic hypersphere or to a non-homogeneous real hypersurface, which is considered as a tube of radius $r=\frac{\pi}{4}$ over a holomorphic curve. In case of \CHc, \M is locally congruent to a Hopf hypersurface with $A\xi=0$ (see \cite{IR2009}) and this completes the proof of Theorem \ref{txiparallel}.

\begin{remark}
In case of real hypersurfaces with constant $\lambda,\nu$  in direction of $\xi$ it can be proved that the eigenvalues of real hypersurfaces of type ($B$) both in \CPc and \CHc satisfies the above. So the *-Ricci tensor of real hypesurfaces of type ($B$) is $\xi$-parallel.
\end{remark}

\scriptsize{\textsc{\hspace{-15pt}G. Kaimakamis, Faculty of Mathematics and Engineering Sciences, Hellenic Military Academy, Vari, Attiki, Greece}\\
\textsc{e-mail}: gmiamis@gmail.com\\

\textsc{\hspace{-15pt}K. Panagiotidou, Faculty of Mathematics and Engineering Sciences, Hellenic Military Academy, Vari, Attiki, Greece}}\\
\textsc{e-mail}: konpanagiotidou@gmail.com\\

\end{document}